\def\temp{1.34}%
\let\tempp=\relax
\expandafter\ifx\csname psboxversion\endcsname\relax
  \message{PSBOX(\temp) loading}%
\else
    \ifdim\temp cm>\psboxversion cm
      \message{PSBOX(\temp) loading}%
    \else
      \message{PSBOX(\psboxversion) is already loaded: I won't load
        PSBOX(\temp)!}%
      \let\temp=\psboxversion
      \let\tempp= 
    \fi
\fi
\tempp
\let\psboxversion=\temp
\catcode`\@=11
%
%
\def\psfortextures{
\def\PSspeci@l##1##2{%
\special{illustration ##1\space scaled ##2}%
}}%
\def\psfordvitops{
\def\PSspeci@l##1##2{%
\special{dvitops: import ##1\space \the\drawingwd \the\drawinght}%
}}%
\def\psfordvips{
\def\PSspeci@l##1##2{%
\d@my=0.1bp \d@mx=\drawingwd \divide\d@mx by\d@my
\includegraphics{##1\space}}}%
\def\psforoztex{
\def\PSspeci@l##1##2{%
\special{##1 \space
      ##2 1000 div dup scale
      \number-\psllx\space \number-\pslly\space translate
}}}%
\def\psfordvitps{
\def\psdimt@n@sp##1{\d@mx=##1\relax\edef\psn@sp{\number\d@mx}}
\def\PSspeci@l##1##2{%
\special{dvitps: Include0 "psfig.psr"}
\psdimt@n@sp{\drawingwd}
\special{dvitps: Literal "\psn@sp\space"}
\psdimt@n@sp{\drawinght}
\special{dvitps: Literal "\psn@sp\space"}
\psdimt@n@sp{\psllx bp}
\special{dvitps: Literal "\psn@sp\space"}
\psdimt@n@sp{\pslly bp}
\special{dvitps: Literal "\psn@sp\space"}
\psdimt@n@sp{\psurx bp}
\special{dvitps: Literal "\psn@sp\space"}
\psdimt@n@sp{\psury bp}
\special{dvitps: Literal "\psn@sp\space startTexFig\space"}
\special{dvitps: Include1 "##1"}
\special{dvitps: Literal "endTexFig\space"}
}}%
\def\psfordvialw{
\def\PSspeci@l##1##2{
\special{language "PostScript",
position = "bottom left",
literal "  \psllx\space \pslly\space translate
  ##2 1000 div dup scale
  -\psllx\space -\pslly\space translate",
include "##1"}
}}%
\def\psforptips{
\def\PSspeci@l##1##2{{
\d@mx=\psurx bp
\advance \d@mx by -\psllx bp
\divide \d@mx by 1000\multiply\d@mx by \xscale
\incm{\d@mx}
\let\tmpx\dimincm
\d@my=\psury bp
\advance \d@my by -\pslly bp
\divide \d@my by 1000\multiply\d@my by \xscale
\incm{\d@my}
\let\tmpy\dimincm
\d@mx=-\psllx bp
\divide \d@mx by 1000\multiply\d@mx by \xscale
\d@my=-\pslly bp
\divide \d@my by 1000\multiply\d@my by \xscale
\at(\d@mx;\d@my){\special{ps:##1 x=\tmpx, y=\tmpy}}
}}}%
\def\psonlyboxes{
\def\PSspeci@l##1##2{%
\at(0cm;0cm){\boxit{\vbox to\drawinght
  {\vss\hbox to\drawingwd{\at(0cm;0cm){\hbox{({\tt##1})}}\hss}}}}
}}%
\def\psloc@lerr#1{%
\let\savedPSspeci@l=\PSspeci@l%
\def\PSspeci@l##1##2{%
\at(0cm;0cm){\boxit{\vbox to\drawinght
  {\vss\hbox to\drawingwd{\at(0cm;0cm){\hbox{({\tt##1}) #1}}\hss}}}}
\let\PSspeci@l=\savedPSspeci@l
}}%
%
%
\newread\pst@mpin
\newdimen\drawinght\newdimen\drawingwd
\newdimen\psxoffset\newdimen\psyoffset
\newbox\drawingBox
\newcount\xscale \newcount\yscale \newdimen\pscm\pscm=1cm
\newdimen\d@mx \newdimen\d@my
\newdimen\pswdincr \newdimen\pshtincr
\let\ps@nnotation=\relax
{\catcode`\|=0 |catcode`|\=12 |catcode`|
|catcode`#=12 |catcode`*=14
|xdef|backslashother{\}*
|xdef|percentother{
|xdef|tildeother{~}*
|xdef|sharpother{#}*
}%
\def\R@moveMeaningHeader#1:->{}%
\def\uncatcode#1{%
\edef#1{\expandafter\R@moveMeaningHeader\meaning#1}}%
\def\execute#1{#1}
\def\psm@keother#1{\catcode`#112\relax}
\def\executeinspecs#1{%
\execute{\begingroup\let\do\psm@keother\dospecials\catcode`\^^M=9#1\endgroup}}%
\def\@mpty{}%
\def\matchexpin#1#2{
  \fi%
  \edef\tmpb{{#2}}%
  \expandafter\makem@tchtmp\tmpb%
  \edef\tmpa{#1}\edef\tmpb{#2}%
  \expandafter\expandafter\expandafter\m@tchtmp\expandafter\tmpa\tmpb\endm@tch%
  \if\match%
}%
\def\matchin#1#2{%
  \fi%
  \makem@tchtmp{#2}%
  \m@tchtmp#1#2\endm@tch%
  \if\match%
}%
\def\makem@tchtmp#1{\def\m@tchtmp##1#1##2\endm@tch{%
  \def\tmpa{##1}\def\tmpb{##2}\let\m@tchtmp=\relax%
  \ifx\tmpb\@mpty\def\match{YN}%
  \else\def\match{YY}\fi%
}}%
\def\incm#1{{\psxoffset=1cm\d@my=#1
 \d@mx=\d@my
  \divide\d@mx by \psxoffset
  \xdef\dimincm{\number\d@mx.}
  \advance\d@my by -\number\d@mx cm
  \multiply\d@my by 100
 \d@mx=\d@my
  \divide\d@mx by \psxoffset
  \edef\dimincm{\dimincm\number\d@mx}
  \advance\d@my by -\number\d@mx cm
  \multiply\d@my by 100
 \d@mx=\d@my
  \divide\d@mx by \psxoffset
  \xdef\dimincm{\dimincm\number\d@mx}
}}%
%
\newif\ifNotB@undingBox
\newhelp\PShelp{Proceed: you'll have a 5cm square blank box instead of
your graphics (Jean Orloff).}%
\def\s@tsize#1 #2 #3 #4\@ndsize{
  \def\psllx{#1}\def\pslly{#2}%
  \def\psurx{#3}\def\psury{#4}
  \ifx\psurx\@mpty\NotB@undingBoxtrue
  \else
    \drawinght=#4bp\advance\drawinght by-#2bp
    \drawingwd=#3bp\advance\drawingwd by-#1bp
  \fi
  }%
\def\sc@nBBline#1:#2\@ndBBline{\edef\p@rameter{#1}\edef\v@lue{#2}}%
\def\g@bblefirstblank#1#2:{\ifx#1 \else#1\fi#2}%
{\catcode`\%=12
\xdef\B@undingBox{
\def\ReadPSize#1{
 \readfilename#1\relax
 \let\PSfilename=\lastreadfilename
 \openin\pst@mpin=#1\relax
 \ifeof\pst@mpin \errhelp=\PShelp
   \errmessage{I haven't found your postscript file (\PSfilename)}%
   \psloc@lerr{was not found}%
   \s@tsize 0 0 142 142\@ndsize
   \closein\pst@mpin
 \else
   \if\matchexpin{\GlobalInputList}{, \lastreadfilename}%
   \else\xdef\GlobalInputList{\GlobalInputList, \lastreadfilename}%
     \immediate\write\psbj@inaux{\lastreadfilename,}%
   \fi%
   \loop
     \executeinspecs{\catcode`\ =10\global\read\pst@mpin to\n@xtline}%
     \ifeof\pst@mpin
       \errhelp=\PShelp
       \errmessage{(\PSfilename) is not an Encapsulated PostScript File:
           I could not find any \B@undingBox: line.}%
       \edef\v@lue{0 0 142 142:}%
       \psloc@lerr{is not an EPSFile}%
       \NotB@undingBoxfalse
     \else
       \expandafter\sc@nBBline\n@xtline:\@ndBBline
       \ifx\p@rameter\B@undingBox\NotB@undingBoxfalse
         \edef\t@mp{%
           \expandafter\g@bblefirstblank\v@lue\space\space\space}%
         \expandafter\s@tsize\t@mp\@ndsize
       \else\NotB@undingBoxtrue
       \fi
     \fi
   \ifNotB@undingBox\repeat
   \closein\pst@mpin
 \fi
\message{#1}%
}%
%
%
\def\psboxto(#1;#2)#3{\vbox{%
   \ReadPSize{#3}%
   \advance\pswdincr by \drawingwd
   \advance\pshtincr by \drawinght
   \divide\pswdincr by 1000
   \divide\pshtincr by 1000
   \d@mx=#1
   \ifdim\d@mx=0pt\xscale=1000
         \else \xscale=\d@mx \divide \xscale by \pswdincr\fi
   \d@my=#2
   \ifdim\d@my=0pt\yscale=1000
         \else \yscale=\d@my \divide \yscale by \pshtincr\fi
   \ifnum\yscale=1000
         \else\ifnum\xscale=1000\xscale=\yscale
                    \else\ifnum\yscale<\xscale\xscale=\yscale\fi
              \fi
   \fi
   \divide\drawingwd by1000 \multiply\drawingwd by\xscale
   \divide\drawinght by1000 \multiply\drawinght by\xscale
   \divide\psxoffset by1000 \multiply\psxoffset by\xscale
   \divide\psyoffset by1000 \multiply\psyoffset by\xscale
   \global\divide\pscm by 1000
   \global\multiply\pscm by\xscale
   \multiply\pswdincr by\xscale \multiply\pshtincr by\xscale
   \ifdim\d@mx=0pt\d@mx=\pswdincr\fi
   \ifdim\d@my=0pt\d@my=\pshtincr\fi
   \message{scaled \the\xscale}%
 \hbox to\d@mx{\hss\vbox to\d@my{\vss
   \global\setbox\drawingBox=\hbox to 0pt{\kern\psxoffset\vbox to 0pt{%
      \kern-\psyoffset
      \PSspeci@l{\PSfilename}{\the\xscale}%
      \vss}\hss\ps@nnotation}%
   \global\wd\drawingBox=\the\pswdincr
   \global\ht\drawingBox=\the\pshtincr
   \global\drawingwd=\pswdincr
   \global\drawinght=\pshtincr
   \baselineskip=0pt
   \copy\drawingBox
 \vss}\hss}%
  \global\psxoffset=0pt
  \global\psyoffset=0pt
  \global\pswdincr=0pt
  \global\pshtincr=0pt 
  \global\pscm=1cm 
}}%
%
%
\def\psboxscaled#1#2{\vbox{%
  \ReadPSize{#2}%
  \xscale=#1
  \message{scaled \the\xscale}%
  \divide\pswdincr by 1000 \multiply\pswdincr by \xscale
  \divide\pshtincr by 1000 \multiply\pshtincr by \xscale
  \divide\psxoffset by1000 \multiply\psxoffset by\xscale
  \divide\psyoffset by1000 \multiply\psyoffset by\xscale
  \divide\drawingwd by1000 \multiply\drawingwd by\xscale
  \divide\drawinght by1000 \multiply\drawinght by\xscale
  \global\divide\pscm by 1000
  \global\multiply\pscm by\xscale
  \global\setbox\drawingBox=\hbox to 0pt{\kern\psxoffset\vbox to 0pt{%
     \kern-\psyoffset
     \PSspeci@l{\PSfilename}{\the\xscale}%
     \vss}\hss\ps@nnotation}%
  \advance\pswdincr by \drawingwd
  \advance\pshtincr by \drawinght
  \global\wd\drawingBox=\the\pswdincr
  \global\ht\drawingBox=\the\pshtincr
  \global\drawingwd=\pswdincr
  \global\drawinght=\pshtincr
  \baselineskip=0pt
  \copy\drawingBox
  \global\psxoffset=0pt
  \global\psyoffset=0pt
  \global\pswdincr=0pt
  \global\pshtincr=0pt 
  \global\pscm=1cm
}}%
%
\def\psbox#1{\psboxscaled{1000}{#1}}%
\newif\ifn@teof\n@teoftrue
\newif\ifc@ntrolline
\newif\ifmatch
\newread\j@insplitin
\newwrite\j@insplitout
\newwrite\psbj@inaux
\immediate\openout\psbj@inaux=psbjoin.aux
\immediate\write\psbj@inaux{\string\joinfiles}%
\immediate\write\psbj@inaux{\jobname,}%
%
%
\def\toother#1{\ifcat\relax#1\else\expandafter%
  \toother@ux\meaning#1\endtoother@ux\fi}%
\def\toother@ux#1 #2#3\endtoother@ux{\def\tmp{#3}%
  \ifx\tmp\@mpty\def\tmp{#2}\let\next=\relax%
  \else\def\next{\toother@ux#2#3\endtoother@ux}\fi%
\next}%
%
%
\let\readfilenamehook=\relax
\def\re@d{\expandafter\re@daux}
\def\re@daux{\futurelet\nextchar\stopre@dtest}%
\def\re@dnext{\xdef\lastreadfilename{\lastreadfilename\nextchar}%
  \afterassignment\re@d\let\nextchar}%
\def\stopre@d{\egroup\readfilenamehook}%
\def\stopre@dtest{%
  \ifcat\nextchar\relax\let\nextread\stopre@d
  \else
    \ifcat\nextchar\space\def\nextread{%
      \afterassignment\stopre@d\chardef\nextchar=`}%
    \else\let\nextread=\re@dnext
      \toother\nextchar
      \edef\nextchar{\tmp}%
    \fi
  \fi\nextread}%
\def\readfilename{\bgroup%
  \let\\=\backslashother \let\%=\percentother \let\~=\tildeother
  \let\#=\sharpother \xdef\lastreadfilename{}%
  \re@d}%
%
%
\xdef\GlobalInputList{\jobname}%
\def\psnewinput{%
  \def\readfilenamehook{
    \if\matchexpin{\GlobalInputList}{, \lastreadfilename}%
    \else\xdef\GlobalInputList{\GlobalInputList, \lastreadfilename}%
      \immediate\write\psbj@inaux{\lastreadfilename,}%
    \fi%
    \ps@ldinput\lastreadfilename\relax%
    \let\readfilenamehook=\relax%
  }\readfilename%
}%
\expandafter\ifx\csname @@input\endcsname\relax    
  \immediate\let\ps@ldinput=\input\def\input{\psnewinput}%
\else
  \immediate\let\ps@ldinput=\@@input
  \def\@@input{\psnewinput}%
\fi%
\def\nowarnopenout{%
 \def\warnopenout##1##2{%
   \readfilename##2\relax
   \message{\lastreadfilename}%
   \immediate\openout##1=\lastreadfilename\relax}}%
\def\warnopenout#1#2{%
 \readfilename#2\relax
 \def\t@mp{TrashMe,psbjoin.aux,psbjoint.tex,}\uncatcode\t@mp
 \if\matchexpin{\t@mp}{\lastreadfilename,}%
 \else
   \immediate\openin\pst@mpin=\lastreadfilename\relax
   \ifeof\pst@mpin
     \else
     \errhelp{If the content of this file is so precious to you, abort (ie
press x or e) and rename it before retrying.}%
     \errmessage{I'm just about to replace your file named \lastreadfilename}%
   \fi
   \immediate\closein\pst@mpin
 \fi
 \message{\lastreadfilename}%
 \immediate\openout#1=\lastreadfilename\relax}%
{\catcode`\%=12\catcode`\*=14
\gdef\splitfile#1{*
 \readfilename#1\relax
 \immediate\openin\j@insplitin=\lastreadfilename\relax
 \ifeof\j@insplitin
   \message{! I couldn't find and split \lastreadfilename!}*
 \else
   \immediate\openout\j@insplitout=TrashMe
   \message{< Splitting \lastreadfilename\space into}*
   \loop
     \ifeof\j@insplitin
       \immediate\closein\j@insplitin\n@teoffalse
     \else
       \n@teoftrue
       \executeinspecs{\global\read\j@insplitin to\spl@tinline\expandafter
         \ch@ckbeginnewfile\spl@tinline
       \ifc@ntrolline
       \else
         \toks0=\expandafter{\spl@tinline}*
         \immediate\write\j@insplitout{\the\toks0}*
       \fi
     \fi
   \ifn@teof\repeat
   \immediate\closeout\j@insplitout
 \fi\message{>}*
}*
\gdef\ch@ckbeginnewfile#1
 \def\t@mp{#1}*
 \ifx\@mpty\t@mp
   \def\t@mp{#3}*
   \ifx\@mpty\t@mp
     \global\c@ntrollinefalse
   \else
     \immediate\closeout\j@insplitout
     \warnopenout\j@insplitout{#2}*
     \global\c@ntrollinetrue
   \fi
 \else
   \global\c@ntrollinefalse
 \fi}*
\gdef\joinfiles#1\into#2{*
 \message{< Joining following files into}*
 \warnopenout\j@insplitout{#2}*
 \message{:}*
 {*
 \edef\w@##1{\immediate\write\j@insplitout{##1}}*
\w@{
\w@{
\w@{
\w@{
\w@{
\w@{
\w@{
\w@{
\w@{
\w@{
\w@{\string\input\space psbox.tex}*
\w@{\string\splitfile{\string\jobname}}*
\w@{\string\let\string\autojoin=\string\relax}*
}*
 \expandafter\tre@tfilelist#1, \endtre@t
 \immediate\closeout\j@insplitout
 \message{>}*
}*
\gdef\tre@tfilelist#1, #2\endtre@t{*
 \readfilename#1\relax
 \ifx\@mpty\lastreadfilename
 \else
   \immediate\openin\j@insplitin=\lastreadfilename\relax
   \ifeof\j@insplitin
     \errmessage{I couldn't find file \lastreadfilename}*
   \else
     \message{\lastreadfilename}*
     \immediate\write\j@insplitout{
     \executeinspecs{\global\read\j@insplitin to\oldj@ininline}*
     \loop
       \ifeof\j@insplitin\immediate\closein\j@insplitin\n@teoffalse
       \else\n@teoftrue
         \executeinspecs{\global\read\j@insplitin to\j@ininline}*
         \toks0=\expandafter{\oldj@ininline}*
         \let\oldj@ininline=\j@ininline
         \immediate\write\j@insplitout{\the\toks0}*
       \fi
     \ifn@teof
     \repeat
   \immediate\closein\j@insplitin
   \fi
   \tre@tfilelist#2, \endtre@t
 \fi}*
}%
\def\autojoin{%
 \immediate\write\psbj@inaux{\string\into{psbjoint.tex}}%
 \immediate\closeout\psbj@inaux
 \expandafter\joinfiles\GlobalInputList\into{psbjoint.tex}%
}%
%
%
%
\def\centinsert#1{\midinsert\line{\hss#1\hss}\endinsert}%
\def\psannotate#1#2{\vbox{%
  \def\ps@nnotation{#2\global\let\ps@nnotation=\relax}#1}}%
\def\pscaption#1#2{\vbox{%
   \setbox\drawingBox=#1
   \copy\drawingBox
   \vskip\baselineskip
   \vbox{\hsize=\wd\drawingBox\setbox0=\hbox{#2}%
     \ifdim\wd0>\hsize
       \noindent\unhbox0\tolerance=5000
    \else\centerline{\box0}%
    \fi
}}}%
%
\def\at(#1;#2)#3{\setbox0=\hbox{#3}\ht0=0pt\dp0=0pt
  \rlap{\kern#1\vbox to0pt{\kern-#2\box0\vss}}}%
%
\newdimen\gridht \newdimen\gridwd
\def\gridfill(#1;#2){%
  \setbox0=\hbox to 1\pscm
  {\vrule height1\pscm width.4pt\leaders\hrule\hfill}%
  \gridht=#1
  \divide\gridht by \ht0
  \multiply\gridht by \ht0
  \gridwd=#2
  \divide\gridwd by \wd0
  \multiply\gridwd by \wd0
  \advance \gridwd by \wd0
  \vbox to \gridht{\leaders\hbox to\gridwd{\leaders\box0\hfill}\vfill}}%
%
\def\fillinggrid{\at(0cm;0cm){\vbox{%
  \gridfill(\drawinght;\drawingwd)}}}%
%
%
\def\textleftof#1:{%
  \setbox1=#1
  \setbox0=\vbox\bgroup
    \advance\hsize by -\wd1 \advance\hsize by -2em}%
\def\textrightof#1:{%
  \setbox0=#1
  \setbox1=\vbox\bgroup
    \advance\hsize by -\wd0 \advance\hsize by -2em}%
\def\endtext{%
  \egroup
  \hbox to \hsize{\valign{\vfil##\vfil\cr%
\box0\cr%
\noalign{\hss}\box1\cr}}}%
%
\def\frameit#1#2#3{\hbox{\vrule width#1\vbox{%
  \hrule height#1\vskip#2\hbox{\hskip#2\vbox{#3}\hskip#2}%
        \vskip#2\hrule height#1}\vrule width#1}}%
\def\boxit#1{\frameit{0.4pt}{0pt}{#1}}%
\catcode`\@=12 
%
 \psfordvips   

\magnification=\magstep1
\tolerance=5000
\parskip 3pt
\font\Bbb=msbm10
\font\Bbbs=msbm8
\textfont12=\Bbb
\scriptfont12=\Bbbs
\font\rmsmall= cmr8
\font\bfone= cmbx10 scaled \magstep1
\let\mcd=\mathchardef
\mcd\Ee="7C45
\mcd\Pe="7C50
\mcd\Re="7C52
\mcd\Ze="7C5A
\def\al{\alpha}
\def\be{\beta}
\def\del{\delta}
\def\ep{\epsilon}
\def\ga{\gamma}
\def\Ga{\Gamma}
\def\la{\lambda}
\def\vphi{\varphi}
\def\bra{\langle}
\def\ket{\rangle}
\def\O{{\cal O}}
\def\nea{\nearrow}
\def\sea{\searrow}
\def\12{{\textstyle{1\over2}}}
\centerline{\bfone Random walk versus random line}
\bigskip
\centerline{Jo\"el De Coninck\footnote{$^{(1)}$}
{\rmsmall Centre de Recherche en Mod\'elisation Mol\'eculaire, Universit\'e de 
Mons-Hainaut, 20 Place du Parc, 7000 Mons, Belgium. 
Email: Joel.De.Coninck@crmm.umh.ac.be}, 
Fran\c cois Dunlop\footnote{$^{(2)}$}
{\rmsmall Laboratoire de Physique Th\'eorique et Mod\'elisation (CNRS - UMR 
8089), Universit\'e de Cergy-Pontoise, 95302 Cergy-Pontoise, France. 
Email: Francois.Dunlop@u-cergy.fr, Thierry.Huillet@u-cergy.fr}, 
Thierry Huillet$^{(2)}$}
\bigskip\noindent{\bf Abstract:} {\rm
We consider random walks $X_n$ in $\Ze_+$, obeying a detailed balance
condition, with a weak drift towards the origin when
$X_n\nea\infty$. We reconsider the 
equivalence in law between a random walk bridge and a 1+1 dimensional
Solid-On-Solid bridge with a corresponding Hamiltonian. 
Phase diagrams are discussed in terms of recurrence versus wetting.
A drift $-\del X_n^{-1}+\O(X_n^{-2})$ of the random walk yields a
Solid-On-Solid potential with an attractive well at the origin and a
repulsive tail ${\del(2+\del)\over8}X_n^{-2}+\O(X_n^{-3})$ at infinity,
showing complete wetting for $\del\le1$ and 
critical partial wetting for $\del>1$.}
 
\medskip\noindent{\rmsmall KEYWORDS: Random walk, recurrence, SOS
model, pinning, wetting}  
 
\bigskip\noindent
{\bf 1. Introduction}
\medskip\noindent
We consider a random walk on $\Ze$ or $\Ze_+$ as defined by transition
probabilities 
$\Pe(X_{n+1}|X_n)$, so that the probability distribution of a random
walk bridge of length $N$ is
$$
\Pe(X_1,\dots,X_{N-1}|X_0=X_N=0)=\prod_{n=0}^{N-1}\Pe(X_{n+1}|X_n)
\Big/\Pe(X_N=0|X_0=0) \eqno{(1.1)} 
$$
We consider a random line making a bridge of length $N$, in the form
of a Solid-On-Solid model, as defined 
by a probability distribution of the form
$$
\Pe(X_1,\dots,X_{N-1}|X_0=X_N=0)=Z_N^{-1}\prod_{n=0}^{N-1}e^{-W(X_n,X_{n+1})}
\prod_{n=1}^{N}e^{-V(X_n)} \eqno{(1.2)}
$$
with $W(X,Y)=W(Y,X)$ for all $X,Y$, and $Z_N$ the
partition function normalising the probability.

We address the question of translating $\Pe(X_{n+1}|X_n)$ into
$W(X,Y)$ and $V(X)$, and conversely, and transferring the information
about transience / null recurrence / positive recurrence of the walk
to complete wetting / partial wetting of the SOS model, and back.
This question is related to the Hamiltonian on
random walk trajectories in Ferrari-Mart\'\i nez [FM].

We assume that the walk obeys the detailed balance condition with
respect to a measure on $\Ze$, not necessarily normalisable, which we
write as $\exp(-U(X))$, so that
$$
\eqalign{
\Pe(X_{n+1}|X_n)&=e^{U(X_n)-U(X_{n+1})}\,\Pe(X_n|X_{n+1})\cr
&=e^{\12 U(X_n)}\Bigl(\Pe(X_{n+1}|X_n)\Pe(X_n|X_{n+1})\Bigr)^\12
e^{-\12 U(X_{n+1})}\cr
&\equiv e^{\12 U(X_n)}e^{-W(X_n,X_{n+1})}e^{-\12 U(X_{n+1})}
} \eqno{(1.3)}$$
which defines $W(X,Y)$ from $\Pe(Y|X)$.
The probability of a random walk bridge may now be written as
$$
\Pe(X_1,\dots,X_{N-1}|X_N=X_0=0)=Z_N^{-1}\prod_{n=0}^{N-1}e^{-W(X_n,X_{n+1})}
\eqno{(1.4)}
$$
with $Z_N=\Pe(X_N=0|X_0=0)$, which is of the form (1.2). The detailed balance
condition was used, but the formula implied by (1.3) for the resulting
SOS interaction $W$ does not require the knowledge of the invariant
measure $\exp(-U(X))$. The interaction $W(X_n,X_{n+1})$
typically contains a part of the form $(V(X_n)+V(X_{n+1}))/2$, which may
be split from $W$.

Conversely, given a SOS probability distribution of the form (1.2),
where we let $W$ absorb $V$ like in (1.4),
we look for a set of random walk probability transitions of the form
$$
\Pe(X_{n+1}|X_n)={e^{-W(X_n,X_{n+1})-\12U(X_{n+1})+\12 U(X_n)}\over Z(X_n)}
\eqno{(1.5)}
$$
These would lead to
$$
\Pe(X_1,\dots,X_{N-1}|X_0=X_N=0)=
\prod_{n=0}^{N-1}{e^{-W(X_n,X_{n+1})}\over Z(X_n)}\eqno{(1.6)}
$$
which agrees with (1.4) only if $Z=\,$const., which requires $\exp(-\12 U)$ 
to be an eigenvector of the symmetric kernel $\exp(-W(X,Y))$:
$$
\sum_{X}e^{-\12U(X)}e^{-W(X,Y)}=\rho\,e^{-\12 U(Y)}\eqno{(1.7)}
$$
The Perron-Frobenius theorem [S] indicates that (1.7) should have a
solution $(\rho,U)$. In any case,
(1.7) is equivalent to $\exp(- U)$ being a left-eigenvector of the
(non-symmetric) kernel (1.5) with $Z=\,$const.:
$$
\sum_{X}e^{-U(X)}e^{-W(X,Y)-\12U(Y)+\12 U(X)}
=\rho\,e^{-U(Y)}\eqno{(1.8)}
$$
Therefore (1.5) with $Z=\,$const. and $U$ obeying (1.7) or (1.8)
is an answer to formulating an SOS random line with probability (1.2),
written as (1.4), in terms of a random walk. 
However, it does require the knowledge of the measure $\exp(-U(X))$,
with respect to which the walk will obey the detailed balance condition.
This is related to the transfer matrix solution of the 1+1 dimensional
SOS models of wetting derived in the early eighties [AD, Bu, C, CW, LH, VL] and
further elaborated with path space limit theorems in the late nineties
[Bo, DGZ, IY, V] and references therein. 
Expressing an SOS bridge in terms of a random walk, asymptotically as
$N\to\infty$, was used also in the proof of the Wulff shape for SOS
models (Theorem 1 in [DDR]). 

In the following sections we consider examples, translating from random walk
to SOS model, when
$$
\Pe(X_{n+1}<X_n|X_n)-\Pe(X_{n+1}>X_n|X_n)\sim{\del\over X_n}
\qquad{\rm as}\qquad X_n\to\infty\eqno{(1.9)}
$$
and discuss recurrence versus wetting. 
Interest into such random walks goes back to Lamperti [L1, L2].
Detailed properties of the
random walk are available [DDH, H] in special instances of (1.9),
yielding the corresponding properties in the corresponding SOS models.
Some of these examples admit constructions for bridges not using the
detailed balance formula. 

\bigskip\noindent
{\bf 2. Bridge with $X_{n+1}-X_n=\pm1$: from random walk to random line}
\medskip\noindent
Let
$$
\vphi:\ \{{\scriptstyle{1\over2},\,{3\over2},\,{5\over2}},\,\dots\}\to\Re
$$
Consider a random walk $X_n$ with state space $\Ze_+=\{0,1,2,\dots\}$,
starting at $X_0=0$, with transition probabilities
$$
\Pe(X_{n+1}|X_n)=
{e^{-(X_{n+1}-X_n)\vphi({X_{n+1}+X_n\over2})}\over
e^{-\vphi(X_n+{1\over2})}+e^{\vphi(X_n-{1\over2})}}\qquad{\rm when}\
X_n\ge1\ {\rm and}\ X_{n+1}=X_n\pm1  \eqno{(2.1)}
$$
and reflection at the origin: $X_{n+1}=1$ whenever $X_n=0$. 
Any random walk with transition probabilities 
$p_x=\Pe(X_{n+1}=x+1|X_n=x)$ and $q_x=1-p_x=\Pe(X_{n+1}=x-1|X_n=x)$ may
be written in the form (2.1): 
take $\vphi(\12)$ arbitrarily, and then solve recursively
$$
\vphi(x+\12)=-\vphi(x-\12)+\ln{q_x\over p_x}\,,\qquad x\ge1 \eqno{(2.2)}
$$ 
From (2.1) we get
$$
\eqalign{
\Pe(X_1,\dots,X_{N-1},X_N=0|X_0=0)&=\prod_{n=0}^{N-1}\Pe(X_{n+1}|X_n)\cr
=\prod_{n=1\atop X_n=0}^Ne^{\vphi({1\over2})}&
\prod_{n=1\atop X_n\ge1}^N
{1\over e^{-\vphi(X_n+{1\over2})}+e^{\vphi(X_n-{1\over2})}}
\prod_{n=0}^{N-1}1_{|X_{n+1}-X_n|=1}\cr
=2^{-N}\prod_{n=1\atop X_n=0}^N2e^{\vphi({1\over2})}&
\prod_{n=1\atop X_n\ge1}^N
{2\over e^{-\vphi(X_n+{1\over2})}+e^{\vphi(X_n-{1\over2})}}
\prod_{n=0}^{N-1}1_{|X_{n+1}-X_n|=1}\cr
=2^{-N}\prod_{n=1}^Ne^{-V(X_n)}&\prod_{n=0}^{N-1}1_{|X_{n+1}-X_n|=1}
} \eqno{(2.3)}
$$
with 
$$
V(X)=-\Bigl(\ln 2+\vphi({1\over2})\Bigr)\,1_{X=0}
+\ln{e^{-\vphi(X+{1\over2})}+e^{\vphi(X-{1\over2})}\over2}\,1_{X\ge1}
\eqno{(2.4)}
$$
The key point in the computation (2.3), instead of using the detailed
balance condition, was the pairing of edge factors, one
factor corresponding to going up the edge and the other factor going
down the edge, leading to the cancellation of factors from the
numerator in (2.1). 
This exact cancellation is restricted to bridges, and requires the
coupling $\vphi$ in (2.1) to be associated with the un-oriented edge 
$\{X_n,X_{n+1}\}$ or to the midpoint $(X_n+X_{n+1})/2$.

\noindent Example (see Fig 1):
$$
\eqalign{
\vphi(x)={\del\over2x}\qquad \Rightarrow\qquad
V(X)&=-(\ln2+\del)\,1_{X=0}
+\ln{e^{-{\del\over2X+1}}+e^{{\del\over2X-1}}\over2}\,1_{X\ge1} \cr
&\sim{\del(2+\del)\over8X^2}\qquad{\rm as}\quad X\to\infty
}\eqno{(2.5)}
$$
Such a potential for $\del>0$, having short range attraction at the wall
and long range repulsion far from the wall, is reminiscent of van der
Waals liquids with a positive Hamaker constant [dG, p846]. The 1+1
dimensional SOS model may be considered a crude effective interface
model where some dimensions and degrees of freedom have been
integrated out in a mean field approximation.
$$
\psboxto(6.6cm;0cm){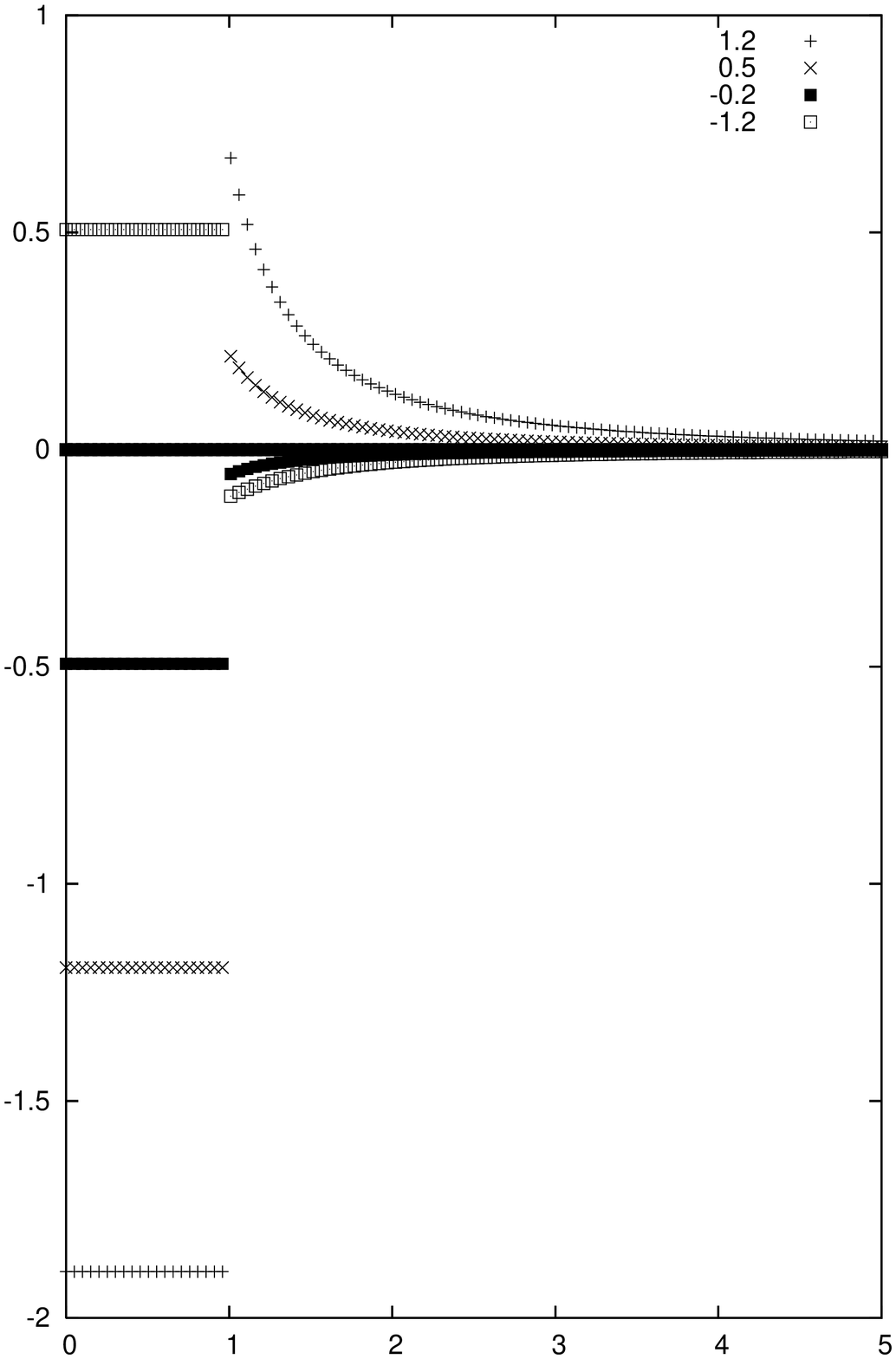}
$$
\vskip-0.2cm\centerline
{\rmsmall Fig. 1: V(X) as (2.4) with $\del=1.2,0.5,-0.2,-1.2$.}

\noindent Example:
$$
\vphi(x)={\del\over2x}+{\ga\over x^2}+\O\Bigl({1\over x^3}\Bigr)
\qquad{\rm as}\quad x\to\infty \eqno{(2.6)}
$$
Such random walks should have a phase diagram (transience / null
recurrence / positive recurrence) independent of $\ga$, and also
independent of the behaviour of $\vphi$ for small $x$. Hence the
corresponding SOS models should have a phase diagram (complete /
partial wetting) independent of $\ga$: partial wetting if and only if 
$\del>1$. However, unlike the square well model in the partial wetting
regime (cf. next section), the height distribution will not decay
exponentially, but as a power law with an exponent depending upon
$\del$ [DDH], hence the term ``critical partial wetting''. 

The behaviour (2.6) implies
$$
V(X)={\del(2+\del)\over8X^2}+\O\Bigl({1\over X^3}\Bigr)
\qquad{\rm as}\quad X\to\infty \eqno{(2.7)}
$$
which indeed is independent of $\ga$.

\bigskip\noindent
{\bf 3. Bridge with $X_{n+1}-X_n=\pm1$: from random line to random walk}
\medskip\noindent
Suppose now that the potential $V(X)$ on $\Ze_+$ is given and
satisfies $V(X)\to0$ as $X\to\infty$. We want to find
$\vphi:\ \{{\scriptstyle{1\over2},\,{3\over2},\,{5\over2}},\,\dots\}\to\Re$
such that (2.4) is satisfied up to a constant $\la$. Let
$$
b_X=e^{V(X)+\la}\,,\qquad a_X=e^{-\vphi(X+\12)} \eqno{(3.1)}
$$
Then (2.4) with $V(X)+\la$ instead of $V(X)$ becomes
$$
\eqalign{
2b_0&=a_0\cr
2b_X&=a_X+a_{X-1}^{-1}\,,\qquad X\ge1
} \eqno{(3.2)}
$$
whose solution is the continued fraction
$$
\eqalign{
a_0&=2b_0\cr
a_1&=2b_1-{1\over2b_0}\cr
&\cdots\cr
a_X&=2b_X-{1\over2b_{X-1}-{1\over2b_{X-2}-{1\over{\cdots\ \cdots\ \cdots\ 
\cdots\ \cdots\ \cdots\ \cdots\ \cdots\atop{\quad\cdots\ \cdots\ \cdots\ 
\cdots\ \cdots\ \cdots\ \cdots
\atop{\qquad{2b_3-{1\over{{2b_2-{1\over{2b_1-{1\over2b_0}}}}}}}}}}}}}
} \eqno{(3.3)}
$$
acceptable only if $a_X>0\ \forall X$. Consistency may be verified
when (2.3)(2.4), converted into (1.4), obeys (1.7), which takes the form
$$
\sum_{X=Y\pm1\atop X\ge0}e^{-\12U(X)-\12V(X)-\12V(Y)}
=2\rho\,e^{-\12 U(Y)}\,,\qquad Y\ge0\eqno{(3.4)}
$$
so that (3.2)(3.3)(3.4) have the solution $\la=\ln\rho$ and
$$
a_X=e^{-\12U(X+1)-\12V(X+1)+\12U(X)+\12V(X)}\,,\qquad X\ge0\eqno{(3.5)}
$$

To conclude this section, we give explicitly the random walks
corresponding to the SOS model with a square well or a double step
potential at the wall: 

\medskip\noindent
$\bullet$ For $V(X)=v_0\,1_{X=0}$, equations (3.1)(3.2) with
$\rho=e^\la$ take the form 
$$
\eqalign{
2b_0&=2\rho e^{v_0}=a_0\cr
2b_X&=2\rho=a_X+a_{X-1}^{-1}\,,\qquad X\ge1
} \eqno{(3.6)}
$$
\item{---}{\sl First ansatz: $\rho=1$}
$$
a_X={(2b_0-1)X+2b_0\over(2b_0-1)X+1}>0\quad \forall X\qquad\Rightarrow\qquad
v_0\ge-\ln2 \eqno{(3.7)}
$$
a transient walk with
$$
\vphi(x)\sim-{1\over x}\qquad{\rm as}\quad x\to\infty \eqno{(3.8)}
$$
compatible with (2.4), $\del=-2$.
\item{---}{\sl Second ansatz: $a_X=a=\,$const.}
$$
\rho={a+a^{-1}\over2}\,,\quad a^{-2}=e^{-v_0}-1>0\qquad\Rightarrow\qquad
v_0<0 \eqno{(3.9)}
$$

\noindent
Both ansatz work when $-\ln2\le v_0<0$, corresponding to transient cases.
The wetting transition is at $v_0=-\ln2$.

\noindent
$\bullet$ For $V(X)=v_0\,1_{X=0}+v_1\,1_{X=1}$, equations (3.1)(3.2) with
$\rho=e^\la$ take the form 
$$
\eqalign{
2b_0&=2\rho e^{v_0}=a_0\cr
2b_1&=2\rho e^{v_1}=a_1+a_0^{-1}\cr
2b_X&=2\rho=a_X+a_{X-1}^{-1}\,,\qquad X\ge2
} \eqno{(3.10)}
$$
\item{---}{\sl First ansatz: $\rho=1$}
$$
\eqalign{
a_0&=2b_0\cr
a_1&=2b_1-{1\over2b_0}>0\cr
a_X&={(a_1-1)X+1\over(a_1-1)X+2-a_1}>0\quad \forall X\ge2
\qquad\Rightarrow\qquad
a_1=2b_1-{1\over2b_0}\ge1
} \eqno{(3.11)}
$$
or
$$
4e^{v_1}\ge2+e^{-v_0} \eqno{(3.12)}
$$
a transient walk with $\vphi(x)\sim-{1\over x}\quad{\rm as}\quad x\to\infty$.
Condition (3.12) coincides with the complete wetting range.

\item{---}{\sl Second ansatz: $a_X=a=\,$const.$\ \forall X\ge1$}
$$
\eqalign{
a_0&=2\rho e^{v_0}\cr
a&=2\rho e^{v_1}-{1\over2\rho e^{v_0}}\cr
\rho&={a+a^{-1}\over2}} \eqno{(3.13)}
$$
Eliminating $\rho$ gives
$$
a^4(e^{v_1}-1)+a^2(2e^{v_1}-e^{-v_0}-1)+e^{v_1}=0 \eqno{(3.14)}
$$
giving a suitable solution for $v_1\le 0$ and any $v_0$ and also for
$$
v_1\ge0\,,\qquad v_0\le 0\,,\qquad v_1\le2\log\cosh{v_0\over2}\eqno{(3.15)}
$$
Whatever $v_0$ and $v_1$, one or the other or both ansatz provides a solution.
There is partial wetting if and only if there is a representation with
$0<a<1$, equivalent to
$$
4e^{v_1}<2+e^{-v_0}\eqno{(3.16)}
$$
where only the second ansatz gives a solution, in fact one solution if
$v_1\le0$ and two solutions if $v_1>0$. 

\vfill\eject
\bigskip
\noindent
{\bf 4. Bridge with $X_{n+1}-X_n\in\{-1,0,+1\}$, Metropolis algorithm}
\medskip\noindent
Let
$$
U:\ \Ze_+\to\Re
$$
Consider a random walk $X_n$ with state space $\Ze_+=\{0,1,2,\dots\}$,
starting at $X_0=0$, with transition probabilities
$$
\eqalign{
\Pe(X_{n+1}|X_n,\,X_n\ge1)=\
&1_{X_{n+1}=X_n\pm1}\,\12e^{-\bigl(U(X_{n+1})-U(X_n)\bigr)_+}\cr
+\,&1_{X_{n+1}=X_n}\,\Bigl[1-\12e^{-\bigl(U(X_n+1)-U(X_n)\bigr)_+}
-\12e^{-\bigl(U(X_n-1)-U(X_n)\bigr)_+}\Bigr]
} \eqno{(4.1)}
$$
and {\sl reflection at the origin}: $X_{n+1}=1$ whenever $X_n=0$. Then
$$
\eqalign{
\Pe(X_1,\dots,X_{N-1},X_N=0|X_0=0)&=\prod_{n=0}^{N-1}\Pe(X_{n+1}|X_n)\cr
=\prod_{n=1\atop X_{n+1}=X_n\ne0}^{N-1}&
\Bigl[1-\12e^{-\bigl(U(X_n+1)-U(X_n)\bigr)_+}
-\12e^{-\bigl(U(X_n-1)-U(X_n)\bigr)_+}\Bigr]\,.\cr
.\,&\prod_{n=1\atop{X_{n+1}=X_n\pm1\atop X_n,X_{n+1}\ne0}}^{N-1}
\12e^{-{|U(X_{n+1})-U(X_n)|\over2}}
\prod_{n=1\atop X_n=0}^N\12e^{-(U(0)-U(1))_+}\cr
=\prod_{n=0\atop X_{n+1}=X_n\ne0}^{N-1}&
\Bigl[1-\12e^{-\bigl(U(X_n+1)-U(X_n)\bigr)_+}
-\12e^{-\bigl(U(X_n-1)-U(X_n)\bigr)_+}\Bigr]\,.\cr
.\,&\prod_{n=0\atop X_{n+1}=X_n\pm1}^{N-1}
\12e^{-{|U(X_{n+1})-U(X_n)|\over2}}
\prod_{n=1\atop X_n=0}^N2e^{(U(1)-U(0))_+}\cr
=2^{-N}\prod_{n=0}^{N-1}&e^{-W(X_n,X_{n+1})}\prod_{n=1}^Ne^{-V(X_n)}
} \eqno{(4.2)}
$$
with 
$$
\eqalign{
V(X)&=-\Bigl(\ln 2+(U(1)-U(0))_+\Bigr)\,1_{X=0}\cr
W(X,X)&=-\ln\Bigl[2-e^{-\bigl(U(X+1)-U(X)\bigr)_+}
-e^{-\bigl(U(X-1)-U(X)\bigr)_+}\Bigr]\cr
W(X,X+1)&=W(X+1,X)={|U(X+1)-U(X)|\over2}
} \eqno{(4.3)}
$$
where $W(X,X)$ is used only with $X\ge1$. The pairing of edge factors
was used, like in Section 2.

Example: $\del\ge0$ and
$$
\eqalign{
U(X)=\del\,\ln(X+1)\qquad \Rightarrow\qquad
V(X)&=-(\ln2+\del\ln2)\,1_{X=0}\cr
W(X,X)&=-\ln\Bigl(1-\Bigl({X+1\over X+2}\Bigr)^\del\Bigr)\cr
W(X+1,X)=W(X,X+1)&={\del\over2}\ln\Bigl({X+2\over X+1}\Bigr)
} \eqno{(4.4)}
$$
\bigskip
Instead of reflection at the origin, let us now choose the full
Metropolis algorithm, including at the wall:
$$
\Pe(X_{n+1}|X_n=0)=\12e^{-(U(1)-U(0))_+}\,1_{X_{n+1}=1}
+(1-\12e^{-(U(1)-U(0))_+})\,1_{X_{n+1}=0} \eqno{(4.5)}
$$
Then
$$
\eqalign{
\Pe(X_1,\dots,X_{N-1},X_N=0|X_0=0)&=\prod_{n=0}^{N-1}\Pe(X_{n+1}|X_n)\cr
=\prod_{n=0\atop X_{n+1}=X_n\ne0}^{N-1}&
\Bigl[1-\12e^{-\bigl(U(X_n+1)-U(X_n)\bigr)_+}
-\12e^{-\bigl(U(X_n-1)-U(X_n)\bigr)_+}\Bigr]\,.\cr
.\,\prod_{n=0\atop X_{n+1}=X_n=0}^{N-1}
&\Bigl[1-\12e^{-\bigl(U(1)-U(0)\bigr)_+}\Bigr]
\prod_{n=0\atop X_{n+1}=X_n\pm1}^{N-1}
\12e^{-{|U(X_{n+1})-U(X_n)|\over2}}\cr
=2^{-N}\prod_{n=0}^{N-1}&e^{-W(X_n,X_{n+1})}
} \eqno{(4.6)}
$$
with 
$$
\eqalign{
W(X,X+1)=W(X+1,X)&={|U(X+1)-U(X)|\over2}\cr
W(X,X)&=-\ln\Bigl[2-e^{-\bigl(U(X+1)-U(X)\bigr)_+}
-e^{-\bigl(U(X-1)-U(X)\bigr)_+}\Bigr]\cr
{\rm except:}\qquad W(0,0)&=-\ln\Bigl[2-e^{-\bigl(U(1)-U(0)\bigr)_+}\Bigr]
} \eqno{(4.7)}
$$
Example:  $\del\ge0$ and
$$
U(X)=\del\,\ln(X+1)\qquad \Rightarrow\qquad
W(0,0)=-\ln\bigl(2-2^{-\del}\bigr) \eqno{(4.8)}
$$
and the other values same as first Metropolis example.

\medskip\noindent
{\bf Remark:} The factor $1/2$ in (4.1) could be replaced by any number
between 0 and $1/2$.

\bigskip\noindent
{\bf 5. Random walk with $X_{n+1}-X_n\in\Ze$, Metropolis algorithm}
\medskip\noindent
Let $\exp(-W_0(X,Y)$ be a symmetric probability kernel in $\Ze\times\Ze$, 
$$
W_0(X,Y)=W_0(Y,X)\,,\qquad\sum_{Y\in\Ze}e^{-W_0(X,Y)}=1 \eqno{(5.1)}
$$
and
$$
U:\ \Ze_+\to\{\Re\cup\{+\infty\}\}\,,\qquad {\rm with:}\quad 
X<0\ \Rightarrow\ U(X)=+\infty \eqno{(5.2)}
$$
Consider a random walk $X_n$ with state space $\Ze_+=\{0,1,2,\dots\}$,
starting at $X_0=0$, with transition probabilities
$$
\eqalign{
\Pe(X_{n+1}|X_n)&=e^{-W_0(X_{n+1},X_n)-\bigl(U(X_{n+1})-U(X_n)\bigr)_+}
\qquad{\rm if}\quad X_{n+1}\ne X_n\cr
\Pe(X_{n+1}=X_n|X_n)&=1-\sum_{Y\ne X_n}\Pe(Y|X_n)
} \eqno{(5.3)}
$$
Then, proceeding as in Section 1, we get (1.4) with
$$
\eqalign{
W(X,Y)=W(Y,X)&=W_0(X,Y)+{|U(Y)-U(X)|\over2}
\qquad{\rm if}\quad Y\ne X\cr
W(X,X)&=-\ln\Bigl(1-\sum_{Y\in\Ze}e^{-W_0(X,Y)-\bigl(U(Y)-U(X)\bigr)_+}\Bigr)
} \eqno{(5.4)}
$$

Example:
$$
W_0(X,Y)=J|X-Y|+{\rm const.}\,,\qquad U(X)=\del\,\ln(X+1) \eqno{(5.5)}
$$
again giving partial wetting if and only if $\del>1$.

Example: $W_0(X,Y)=\ln2$ if $|X-Y|=1$ and $+\infty$ otherwise. This is
equivalent to (4.1)(4.5-7).

For the random walk with $X_{n+1}-X_n\in\Ze$, edges up and down cannot
be paired exactly as in Sections 2 and 4. Approximate pairing would leave
a remainder of order $X_n^{-2}$, which one might argue to be
``irrelevant''.

\bigskip\noindent
{\bf Acknowledgments}: F. D. and T. H. acknowledge support and kind
hospitality from Universit\'e de Mons-Hainaut and CRMM.

\baselineskip 0pt
\bigskip\noindent
{\bf References}

\noindent\item{[AD]} D.B. Abraham, J. De Coninck: {\sl  Description of
phases in a film-thickening transition}, J. Phys. {\bf A16}, L333--337 (1983).

\noindent\item{[Bo]} E. Bolthausen: {\sl  Localization-delocalization phenomena
for random interfaces}, Proceedings of the ICM, Beijing 2002, vol. 3,
pp 25--40.

\noindent\item{[Bu]} T.W. Burkhardt: {\sl Localization--delocalization
transition in a solid-on-solid model with a pinning potential}, J. Phys. 
{\bf A14}, L63--L68 (1981).

\noindent\item{[C]} J.T. Chalker: {\sl  The pinning of a domain wall
by weakened bonds in two dimensions}, J. Phys. {\bf A14}, 2431--2440 (1981).

\noindent\item{[CW]} S.T. Chui, J.D. Weeks: {\sl  Pinning and roughening of
one-dimensional models of interfaces and steps},
Phys. Rev. {\bf B23}, 2438--2441 (1981).

\noindent\item{[dG]} P.G. de Gennes: {\sl  Wetting: statics and dynamics},
Rev. Mod. Phys. {\bf 57}, 827--863 (1985).

\noindent\item{[DDH]} J. De Coninck, F. Dunlop, T. Huillet:
    {\sl Random walk weakly attracted to a wall},  
     J. Stat. Phys. {\bf 133}, 271--280 (2008).

\item{[DDR]} J. De Coninck, F. Dunlop, V. Rivasseau : {\sl On the Microscopic 
Validity of the Wulff Construction and of the Generalized Young Equation}, 
Commun. Math. Phys {\bf 121}, 401--419 (1989).

\noindent\item{[DGZ]} J-D Deuschel, G. Giacomin, L. Zambotti: {\sl
Scaling limits of equilibrium wetting models in (1+1)-dimension},
Probab. Theory Relat. Fields {\bf 132}, 471--500 (2005).

\noindent\item{[FM]} P.A. Ferrari, S. Mart\'\i nez:
    {\sl Hamiltonians on random walk trajectories},
Stochastic Process. Appl.  {\bf 78}, 47--68 (1998).

\noindent\item{[H]} T. Huillet: {\sl 
Random walk with long-range interaction with a
barrier and its dual: Exact results}, Preprint hal-00370353. 

\noindent\item{[IY]} Y. Isozaki, N. Yoshida: {\sl 
Weakly pinned random walk on the wall: pathwise descriptions of the
phase transition}, 
Stochastic Process. Appl. {\bf 96},  261--284 (2001). 

\noindent\item{[L1]} J. Lamperti: {\sl
Criteria for the recurrence or transience of stochastic process I.},
J. Math. Anal. Appl. {\bf 1}, 314--330 (1960).

\noindent\item{[L2]} J. Lamperti: {\sl
Criteria for stochastic processes. II. Passage-time moments},
J. Math. Anal. Appl. {\bf 7}, 127--145 (1963).

\noindent\item{[LH]} J.M.J. van Leeuwen, H.J. Hilhorst: {\sl
Pinning of rough interface by an external potential}, 
Physica A {\bf 107}, 319--329 (1981).

\noindent\item{[S]} E. Seneta: {\sl Nonnegative matrices and Markov chains},
Springer Series in Statistics, Springer-Verlag, New York, 1981.

\noindent\item{[VL]} M. Vallade, J. Lajzerowicz: {\sl 
Transition rugueuse pour une singularit\'e lin\'eaire dans un espace
\`a deux ou trois dimensions},
J. physique {\bf 42}, 1505--1514 (1981).

\noindent\item{[V]} Y. Velenik: {\sl Localization and
delocalization of random interfaces}, Probab. Surv. {\bf 3}, 112--169 (2006).

\bye